\documentclass[letterpaper,12pt]{amsart}

\title{Rigidity theory in statistical inference}

\author{Daniel Irving Bernstein}
\address{Tulane University, Department of Mathematics}
\email{dbernstein1@tulane.edu}

\usepackage{latexsym,array,delarray,amsthm,amssymb,epsfig,mathtools,optidef}
\usepackage[numbers]{natbib}
\usepackage{tikz,url}
\tikzstyle{vertex}=[circle, draw, inner sep=0pt,minimum size=6pt, fill=black]
\newcommand{\vertex}{\node[vertex]}
\usetikzlibrary{calc}
\usepackage{verbatim}
\usepackage{graphicx}
\usepackage{centernot}

\theoremstyle{plain}
\newtheorem{thm}{Theorem}

\newtheorem{conj}[thm]{Conjecture}
\newtheorem*{thm*}{Theorem}
\newtheorem*{lemma*}{Lemma}
\newtheorem*{prop*}{Proposition}
\newtheorem*{cor*}{Corollary}
\newtheorem*{conj*}{Conjecture}

\theoremstyle{definition}
\newtheorem{defn}[thm]{Definition}
\newtheorem*{defn*}{Definition}

\DeclareMathOperator{\tr}{Trace}

\begin{document}

\begin{abstract}
    In this expository article, we summarize what is known about maximum likelihood thresholds of Gaussian models, paying special attention to connections with rigidity theory.
\end{abstract}
\maketitle

\section{Introduction}

A \emph{statistical model} is a collection of probability distributions,
and \emph{model fitting} is the process of choosing the probability distribution within a statistical model
that best describes a particular dataset, according to some criteria.
One of the most fundamental families of statistical models are the \emph{unconstrained $n$-variate Gaussian models},
which are the families of probability distributions whose density functions can be expressed as follows
for some symmetric positive definite $n\times n$ matrix $\Sigma$, called the \emph{covariance matrix} and some $\mu \in \mathbb{R}^n$, called the \emph{mean vector}
\[
    f_{\mu,\Sigma}({\bf x}) := \frac{\exp\left(-\frac{1}{2}({\bf x} -\mu)^T\Sigma^{-1}({\bf x}-\mu)\right)}{\sqrt{(2\pi)^d\det(\Sigma)}}.
\]
In many applications,
one has a dataset ${\bf x}_1,\dots,{\bf x}_d$ consisting of elements of $\mathbb{R}^n$ which are believed to have been independently generated according to some fixed multivariate Gaussian
and the goal is to estimate the parameters $\mu$ and $\Sigma$ which define that distribution.
By treating this dataset as fixed and viewing the entries of $\mu$ and $\Sigma$ as variables,
one can proceed by solving the following optimization problem
\begin{maxi}|l|
    {\mu,\Sigma}{\prod_{i = 1}^d f_{\mu,\Sigma}({\bf x}_i)}
    {}{}
    \addConstraint{\mu \in \mathbb{R}^n {\rm \ and \ } \Sigma {\rm \ positive \ definite}}{}\label{eq: maximum likelihood unconstrained Gaussian}.
\end{maxi}
When \eqref{eq: maximum likelihood unconstrained Gaussian} has a solution, the optimal values $\hat \mu$ and $\hat \Sigma$ for the parameters $\mu$ and $\Sigma$ are the sample mean and covariance, i.e.
\[
    \hat \mu = \frac{1}{d}\sum_{i=1}^d {\bf x}_i \qquad \qquad \hat \Sigma = \sum_{i=1}^d ({\bf x}_i-\hat \mu)({\bf x}_i - \hat\mu)^T.
\]
This solves \eqref{eq: maximum likelihood unconstrained Gaussian} precisely when $\hat \Sigma$ is invertible.

Many statistical models that are used in practice come from adding constraints to an unconstrained Gaussian model.
In particular, for every set $\mathcal{N}$ of $n\times n$ symmetric positive definite matrices, one has the statistical model
whose density functions are the following
\[
    \{f_{\mu,\Sigma} : \mu \in \mathbb{R}^n {\rm \ and \ } \Sigma \in \mathcal{N}\}.
\]
Fitting such a model to a dataset ${\bf x}_1,\dots,{\bf x}_d$ can be done by solving the natural generalization of~\eqref{eq: maximum likelihood unconstrained Gaussian}, namely the following
\begin{maxi}|l|
    {\mu,\Sigma}{\prod_{i = 1}^d f_{\mu,\Sigma}({\bf x}_i)}
    {}{}
    \addConstraint{\mu \in \mathbb{R}^n {\rm \ and \ } \Sigma \in \mathcal{N}}{}\label{eq: constrained maximum likelihood Gaussian}.
\end{maxi}
The mean vector $\mu$ plays no meaningful role in the theoretical questions we are interested in.
By replacing every dataset ${\bf x}_1,\dots,{\bf x}_d$ with its centering around the sample mean $\hat \mu$,
i.e.~with ${\bf x}_1 - \hat \mu,\dots,{\bf x}_d - \hat \mu$,
we can without loss of generality assume that $\mu = 0$.
By doing this and applying several operations to the objective function in~\eqref{eq: constrained maximum likelihood Gaussian} that preserve optimality (taking logarithms, dropping constant terms, etc.), we can transform~\eqref{eq: constrained maximum likelihood Gaussian} into the following equivalent optimization problem where $X$ denotes the $n\times d$ matrix whose columns are the datapoints ${\bf x}_1,\dots,{\bf x}_d$
\begin{mini}|l|
    {\Sigma}{\tr(X^TX \Sigma^{-1}) - \log\det(\Sigma^{-1})}
    {}{}
    \addConstraint{\Sigma \in \mathcal{N}}{}
\end{mini}
We can also shift perspective and view the entries of the \emph{inverse} covariance matrix as our variables.
In particular, we can set
\[
    K := \Sigma^{-1} \qquad {\rm and} \qquad \mathcal{M} := \{A^{-1} : A \in \mathcal{N}\}.
\]
Also notice that $X$ only enters the picture as $X^TX$.
Putting this all together we have an optimization problem of the following form, where $S$ is an arbitrary
$n\times n$ symmetric positive semidefinite matrix of rank at most $d$
and $\mathcal{M}$ is an arbitrary set of $n\times n$ symmetric positive definite matrices
\begin{mini}|l|
    {K}{\tr(SK) - \log\det(K)}
    {}{}
    \addConstraint{K \in \mathcal{M}}{}\label{eq: maximum likelihood final form}.
\end{mini}
This motivates the following definition.
\begin{defn}
    Let $\mathcal{M}$ be a set of $n\times n$ symmetric positive definite matrices.
    The smallest $d$ such that the optimization problem~\eqref{eq: maximum likelihood final form}
    has a solution for almost all $n\times n$ symmetric positive semidefinite $S$ of rank $d$
    is called the \emph{maximum likelihood threshold (MLT) of $\mathcal{M}$}.
\end{defn}
In this article, we will discuss what is known about maximum likelihood thresholds of various families of statistical models $\mathcal{M}$
and how rigidity theory is relevant.

\subsection{Notation}
We use the notation $\mathcal{S}^n,\mathcal{S}^n_+,\mathcal{S}^n_+(d),\mathcal{S}^n_{++}$ to respectively denote the sets of $n\times n$
the sets of symmetric matrices, positive semidefinite matrices, positive semidefinite matrices of rank $d$, and positive definite matrices.
The set $\{1,\dots,n\}$ will be denoted $[n]$.
Given a graph $G$, the set of edges of $G$ will be denoted $E(G)$.

\section{Gaussian graphical models}
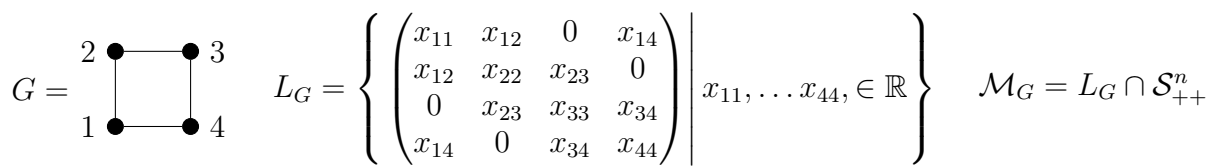
\begin{figure}
    \begin{tikzpicture}
        \vertex (a) at (0,0)[label=left:$1$]{};
        \vertex (b) at (0,1)[label=left:$2$]{};
        \vertex (c) at (1,1)[label=right:$3$]{};
        \vertex (d) at (1,0)[label=right:$4$]{};
        \path
            (a) edge (b)
            (b) edge (c)
            (c) edge (d)
            (d) edge (a)
        ;
        \node at (-1,1/2){$G=$};
        \node at (6.5,1/2){$L_G = \left\{\left.\begin{pmatrix}
        x_{11} & x_{12} & 0 & x_{14} \\
        x_{12} & x_{22} & x_{23} & 0 \\
        0 & x_{23} & x_{33} & x_{34} \\
        x_{14} & 0 & x_{34} & x_{44}
        \end{pmatrix} \right\vert x_{11},\dots x_{44}, \in \mathbb{R} \right\}$};
        \node at (13,1/2){$\mathcal{M}_G = L_G \cap \mathcal{S}^n_{++}$};
    \end{tikzpicture}
    \caption{A graph and its corresponding Gaussian graphical model}\label{fig: graphical model example}
\end{figure}

Let $G$ be a graph on vertex set $[n]$ and define
\[
    \mathcal{M}_G := \{A \in \mathcal{S}^n_{++}: A_{ij} = 0 {\rm \ if \ } i\neq j { \rm \ and \  } ij \notin E(G)\}.
\]
Given a graph $G$,
the corresponding \emph{Gaussian graphical model} is the family of multivariate Gaussian distributions whose inverse covariance lies in $\mathcal{M}_G$.
See Figure~\ref{fig: graphical model example} for an example.
These models were originally introduced and studied by Dempster in 1972 in~\cite{dempster1972covariance}.
A simple consequence of Dempster's result is the following alternative characterization of the MLT of a model $\mathcal{M}_G$.

\begin{thm}[\cite{dempster1972covariance}]\label{thm:dempster}
    Let $G$ be a graph on vertex set $[n]$.
    Then the MLT of $\mathcal{M}_G$ is the minimum $d$ such that for almost every $S \in \mathcal{S}^n_+(d)$,
    there exists $T \in \mathcal{S}^n_{++}$ such that $T_{ii} = S_{ii}$ for all $i \in [n]$ and $T_{ij} = S_{ij}$ for all $ij \in E$.
\end{thm}

Computing or bounding the maximum likelihood threshold of $\mathcal{M}_G$ in terms of the combinatorics of $G$
is surprisingly difficult.
The first result in this direction was proven in 1993 by Buhl in~\cite{buhl1993existence}, which we state below
after recalling the necessary graph theoretic definitions.
Recall that a graph is \emph{chordal} if it has no induced cycles of length four or greater,
and that a \emph{chordal cover} of a graph $G$ is a chordal graph containing $G$ as a subgraph.
The \emph{treewidth} of a graph $G$ is one less than the minimum clique number of a chordal cover of $G$.
The clique number of a chordal graph is therefore one plus its treewidth,
and thus a corollary of the following result is that the MLT of a chordal graph is its clique number.

\begin{thm}[{\cite{buhl1993existence}}]\label{thm: buhl}
    Let $G$ be a graph with clique number $\omega$ and treewidth $\tau$.
    Then the MLT of $G$ lies between $\omega$ and $\tau + 1$.
\end{thm}

In 2012, Uhler showed~\cite{uhler2012geometry} how one can use Gr\"obner bases to bound the MLT of $\mathcal{M}_G$ from above.
Gross and Sullivant later showed~\cite{gross2018maximum} how this upper bound can be understood in terms of rigidity properties of~$G$.
We now describe this after introducing the necessary rigidity-theoretic language.

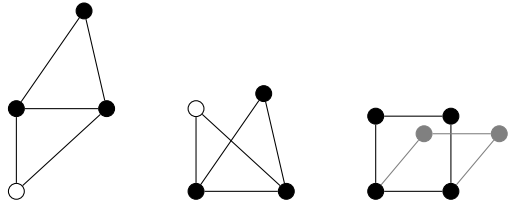
\begin{figure}
    \begin{tikzpicture}
        \vertex (1) at (0,0){};
        \vertex (2) at (1.2,0){};
        \vertex (3) at (.9,1.3){};
        \vertex[fill=white] (4) at (0,-1.1){};
        \path
            (1) edge (2) edge (3) edge (4)
            (2) edge (3) edge (4)
        ;
    \end{tikzpicture}
    \qquad
    \begin{tikzpicture}
        \vertex (1) at (0,0){};
        \vertex (2) at (1.2,0){};
        \vertex (3) at (.9,1.3){};
        \vertex[fill=white] (4) at (0,1.1){};
        \path
            (1) edge (2) edge (3) edge (4)
            (2) edge (3) edge (4)
        ;
    \end{tikzpicture}
    \qquad
    \begin{tikzpicture}
        \vertex (1) at (0,0){};
        \vertex (2) at (1,0){};
        \vertex (3) at (1,1){};
        \vertex (4) at (0,1){};
        \vertex[color=gray] (5) at ({cos(50)},{sin(50)}){};
        \vertex[color=gray] (6) at ({cos(50)+1},{sin(50)}){};
        \path
            (1) edge (2) edge (4)
            (2) edge (3)
            (3) edge (4)
        ;
        \color{gray}
        \path
            (1) edge (5)
            (2) edge (6)
            (5) edge (6)
        ;
    \end{tikzpicture}
    \caption{
        The two frameworks on the left are equivalent but not congruent.
        All edge lengths are the same, but the distance between the two non-adjacent vertices is different.
        One can visualize moving between the frameworks by folding the white vertex out of the plane of the screen.
        Without leaving two dimensions, there is no way to continuously deform one of the frameworks into another equivalent framework,
        so both are locally rigid.
        The framework on the right is not locally rigid since it can be continuously deformed in the plane without changing any edge lengths, as shown.
    }\label{fig:rigidity definitions}
\end{figure}

A \emph{$d$-dimensional framework} consists of a graph $G$ on vertex set $[n]$
and a function $p: [n]\rightarrow \mathbb{R}^d$.
A $d$-dimensional framework $(G,p)$ should be conceptualized as a physical construction of the graph $G$ in $d$-dimensional space
by placing vertex $i$ at point $p(i)$.
Two frameworks $(G,p)$ and $(G,q)$ on the same graph are \emph{equivalent} if $\|p(i)-p(j)\| = \|q(i)-q(j)\|$ for every edge $ij$ of $G$
and moreover \emph{congruent} if $\|p(i)-p(j)\| = \|q(i)-q(j)\|$ holds for every pair $i,j \in [n]$.
The left two frameworks in Figure~\ref{fig:rigidity definitions} are equivalent but not congruent.

A $d$-dimensional framework $(G,p)$ is \emph{locally rigid} if there exists $\varepsilon > 0$ such that whenever
$q: [n]\rightarrow \mathbb{R}^d$ satisfies $\|p(i)-q(i)\| \le \varepsilon$ for all $i \in [n]$,
$(G,p)$ and $(G,q)$ are equivalent.
The frameworks on the left in Figure~\ref{fig:rigidity definitions} are both locally rigid.
The framework on the right is not, and this can be seen by noting that the framework can shear as shown.
For each graph $G$ and each positive integer $d$, there exists a subset $S$ of functions $p: [n]\rightarrow \mathbb{R}^d$ of Lebesge measure zero
such that if $p: [n]\rightarrow \mathbb{R}^d$ is \emph{not} in $S$,
then either $(G,p)$ is always locally rigid, or $(G,p)$ is always not locally rigid.
In the former case, we say that $G$ is \emph{generically locally rigid}.
The four-cycle is not generically locally rigid, but adding an edge produces a graph that is.

The (edge sets of the) graphs that are generically locally rigid in $\mathbb{R}^d$ are the spanning sets of a matroid whose ground set is the edge set of the complete graph.
This matroid is called the \emph{$d$-dimensional rigidity matroid}.
Graphs whose edge sets are independent in this matroid are called \emph{$d$-independent}.
We can now state Gross and Sullivant's rigidity-theoretic characterization of Uhler's bound.

\begin{thm}[{\cite{gross2018maximum}}]\label{thm: upper bound rigidity}
    Let $G$ be a graph on vertex set $[n]$.
    If $G$ is $d$-independent,
    then the maximum likelihood threshold of $\mathcal{M}_G$ is at most $d+1$.
\end{thm}

Theorem~\ref{thm: upper bound rigidity} begs the question of whether the maximum likelihood threshold of a graph $G$ is simply the minimum $d+1$ such that $G$ is $d$-independent.
Blekherman and Sinn answered this question in the negative in~\cite{blekherman2019maximum}.
In particular, they showed that the minimum $d$ for which $K_{5,5}$ is $d$-independent is $4$,
but the MLT of this graph is $4$, not $5$ as an affirmative answer would imply.

One can still ask if the MLT of $\mathcal{M}_G$ can be understood in terms of
rigidity-theoretic properties of $G$.
Dewar, Gortler, Nixon, Theran, Sitharam, and I provided such an understanding in~\cite{bernstein2024maximum}, described in the next theorem.
When we say ``every generic framework,'' we mean that there is a set of functions
$[p]\rightarrow \mathbb{R}^{d-1}$ of Lebesgue measure zero for which the theorem may fail.

\begin{thm}[{\cite[Theorem 1.16]{bernstein2024maximum}}]\label{thm: rigidity classification of MLT}
    Let $G$ be a graph on vertex set $[n]$.
    Then the MLT of $\mathcal{M}_G$ is the minimum dimension $d$ such that every generic $d-1$-dimensional framework $(G,p)$ is equivalent to an $(n-1)$-dimensional framework $(G,q)$
    such that $\{q(i) : i =1,\dots,n\}$ affinely spans $\mathbb{R}^{n-1}$.
\end{thm}

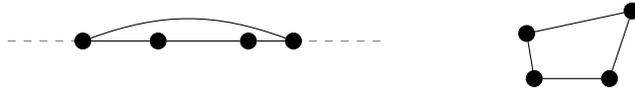
\begin{figure}
    \begin{tikzpicture}
        \draw[gray,dashed] (-1,0.5)--(4,0.5);
        \vertex (1) at (0,0.5){};
        \vertex (2) at (1,0.5){};
        \vertex (3) at (2.2,0.5){};
        \vertex (4) at (2.8,0.5){};
        \path
            (1) edge (2)
            (2) edge (3)
            (3) edge (4)
            (1) edge[bend left = 20] (4)
        ;
        \vertex (1) at (0+6,0){};
        \vertex (2) at (1+6,0){};
        \vertex (3) at (1.3+6,0.9){};
        \vertex (4) at (-0.1+6,0.6){};
        \path
            (1) edge (2) edge (4)
            (2) edge (3)
            (3) edge (4)
        ;   
    \end{tikzpicture}
    \caption{Generic frameworks on the four-cycle in one and two dimensions.
    In the framework on the right, the longest edge length is the sum of the shorter three.
    Any framework satisfying this property must have all vertices lying on a line.
    In particular, the affine span of the vertices of any equivalent framework is one-dimensional.
    Moreover, this property is robust with respect to perturbing the framework (while staying on the line).
    Thus one cannot simply dismiss this obstacle by invoking genericity and so Theorem~\ref{thm: rigidity classification of MLT} implies that the MLT of the four cycle is strictly greater than $2$.
    The two-dimensional framework on the right is equivalent to any three-dimensional framework obtained by folding the framework along a diagonal out of the plane of the page (or computer screen).
    In fact, any generic two-dimensional framework on this graph will satisfy this property.
    Theorem~\ref{thm: rigidity classification of MLT} then implies that the MLT of the four cycle
    is at most 3, and thus exactly 3 in light of the one-dimensional frameworks.}\label{fig: MLT theorem illustration}
\end{figure}

See Figure~\ref{fig: MLT theorem illustration} for an example illustrating Theorem~\ref{thm: rigidity classification of MLT}.
While a completely combinatorial description of the maximum likelihood threshold of Gaussian graphical models remains elusive,
Theorem~\ref{thm: rigidity classification of MLT} is useful for deriving combinatorial classifications of restricted classes of Gaussian graphical models.
Examples of such results can be found in the original paper~\cite{bernstein2024maximum} and the follow-up~\cite{bernstein2024computing}.
The following is a summary of the instances where the upper bound from Theorem~\ref{thm: upper bound rigidity} is known to be sharp.

\begin{thm}{\cite{bernstein2024computing,bernstein2024maximum}}
    Let $G$ be a graph on vertex set $[n]$.
    Then the MLT of $\mathcal{M}_G$ is equal to the smallest $d$ such that $G$ is $d$-independent in each of the following cases:
    \begin{enumerate}
        \item $G$ is $3$-independent
        \item $n \le 9$
        \item $G$ has at most $24$ edges
        \item the complement graph of $G$ has at most $5$ edges
        \item $G$ is connected with minimum degree at most $4$ and maximum degree at most $5$.
    \end{enumerate}
\end{thm}

\section{Random graphs}
Most of the work done on maximum likelihood thresholds of Gaussian graphical models has focused on computing and bounding maximum likelihood thresholds of specific graphs.
The next step is to turn this into a more general understanding of what MLTs to expect from graphs that appear in statistical applications.
Such graphs are often computed from a dataset according to a model selection algorithm and therefore one can view them as random graphs
that depend on whatever random process is generating the data.
This motivates questions about the MLTs of random graphs.
We now discuss two different studies on the MLTs of random graphs,
one of which uses the Erd\H{o}s-Reyni random graph model, which is less realistic but easier to analyze,
and the other of which uses the graphical lasso model selection algorithm applied to datasets generated according to a Gaussian distribution.

\subsection{Erd\H{o}s-Reyni random graphs}
Recall that for an integer $n$ and probability $p$,
the \emph{Erd\H{o}s-Reyni random graph model} $G(n,p)$ is the random process that generates a graph $G$ on $n$ vertices
by including each of the $\binom{n}{2}$ possible edges with probability $p$.
We will consider Erd\H{o}s-Reyni random graph models of the form $G(n,c/n)$ for a fixed positive real number $c$
and consider maximum likelihood thresholds of such graphs as $n\rightarrow \infty$.
In this setting, $c$ is asymptotically the expected average degree of a vertex of a graph generated by $G(n,c/n)$.
The following conjecture is supported by computational experiments presented in~\cite{bernstein2024maximum}.
It says that for Erd\H{o}s-Reyni random graphs, one should expect the upper bound given in Theorem~\ref{thm: upper bound rigidity} to be sharp.

\begin{conj}[\cite{bernstein2024maximum}]\label{conj: upper bound sharp for erdos reyni}
    Let $G$ be a random graph generated according to the model $G(n,c/n)$
    and let $d$ be minimum such that $G$ is $d$-independent.
    Then the MLT of $G$ is $d+1$ with high probability.
\end{conj}

\subsection{Graphical lasso}
If $d < n$, then given $n$ datapoints in $\mathbb{R}^d$,
the problem~\eqref{eq: maximum likelihood final form} will not converge in the unconstrained case (i.e.~when $\mathcal{M} = \mathcal{S}^n_{++}$).
As described in~\cite{friedman2008sparse}, one can add a term to the objective function that penalizes dense inverse covariance matrices,
to obtain a related optimization problem that always converges~\cite{ravikumar2011high-dimensional},
and favors sparse inverse covariance matrices.
In particular if $S$ is the $d\times d$ sample covariance matrix obtained from $n$ centered datapoints
and $\alpha > 0$,
the \emph{graphical lasso optimization problem with parameter $\alpha$} is the following
\begin{mini}|l|
    {K}{\tr(S K) - \log\det(K) + \alpha \sum_{i\neq j} |X_{ij}|}
    {}{}
    \addConstraint{K \in \mathcal{S}^{}n_{++}}{}\label{eq: graphical lasso}.
\end{mini}
The optimum of this is the maximum \emph{a posteriori} estimator
of a certain sparsity model on the inverse covariance matrix~\cite{wang2012bayesian},
though it is known to be statistically inconsistent in certain instances~\cite{heinavaara2016inconsistency}.
Letting $\hat K$ denote the optimum, $\hat K$ lies in the graphical model $\mathcal{M}_{G_{\hat K}}$
where $G_{\hat K}$ is the graph on vertex set $\{1,\dots n\}$ that has an edge from $i$ to $j$ whenever $\hat K_{ij} \neq 0$.
With this in mind, one can view graphical lasso as a model selection algorithm which selects the graphical model $\mathcal{M}_{G_{\hat K}}$.
One can thus ask: when graphical lasso is run on a dataset of size $n$,
what is the probability that the MLT of $G_{\hat{K}}$ is less than or equal to $n$?
In other words, how often does graphical lasso select a model that can be fit via maximum likelihood using the same dataset?
Computational experiments investigating this question can be found in~\cite{bernstein2023maximum}.

\section{Other classes of Gaussian models}
\subsection{Linear concentration models}
Gaussian graphical models are families of $n$-variate Gaussians obtained by constraining particular
off-diagonal entries of the inverse covariance matrix to be zero.
We can generalize this by allowing arbitrary linear constraints on the entries of the inverse covariance matrix.
More formally speaking, let $L \subseteq \mathcal{S}^n$ be a linear subspace
such that $L \cap \mathcal{S}^n_{++} \neq \emptyset$. Then define
\[
    \mathcal{M}_L := L \cap \mathcal{S}^n_{++}.
\]
The \emph{linear concentration model} corresponding to $L$ is the family of multivariate Gaussian distributions whose inverse covariance matrix lies in $\mathcal{M}_L$.
Gaussian graphical models are linear covariance models.
Indeed, if $G$ is a graph on vertex set $[n]$ and $L$ is the linear subspace of $\mathcal{S}^n$ defined
by setting the $ij$ entry to $0$ for each \emph{non}-edge $ij$ of $E$,
then $\mathcal{M}_G = \mathcal{M}_L$.

Remarkably, for most linear concentration models, the maximum likelihood threshold depends only on the dimension of the linear space defining it.
More precisely, for each integer $k \ge 0$ there exists an integer $m(k)$ such that if
if $L$ is randomly sampled according to some ``reasonable'' probability measure on the set of $k$-dimensional linear subspaces of $\mathcal{S}^n$ that contain a positive definite matrix,
then the MLT of $\mathcal{M}_L$ is $m(k)$ with probability $1$.
In particular, there is the following result of Gortler, Theran, and myself.

\begin{thm}[{\cite{bernstein2023maximum}}]\label{thm: linear concentration MLTs}
    Let $L \subseteq \mathcal{S}^n$ be a linear subspace obtained by taking the span
    of $m \le \binom{n+1}{2}$ positive definite matrices sampled independently according to a
    ``reasonable''\footnote{specifically, one that is mutually absolutely continuous with respect to Lebesgue measure}
    probability distribution.
    Then, with probability $1$, the maximum likelihood threshold of $\mathcal{M}_L$
    is the minimum $d$ such that
    \[
        m \le nd - \binom{d}{2}.
    \]
\end{thm}
Given the semi-algebraic nature of the problem (as opposed to merely algebraic), proving Theorem~\ref{thm: linear concentration MLTs} is more involved and topological than an algebraic geometer might expect.

\begin{figure}
    \begin{tikzpicture}
        \vertex (a) at (0,0)[label=left:$1$]{};
        \vertex (b) at (0,1)[label=left:$2$]{};
        \vertex (c) at (1,1)[fill=white,label=right:$3$]{};
        \vertex (d) at (1,0)[fill=white,label=right:$4$]{};
        \path
            (a) edge[style=dashed] (b)
            (b) edge[style=dashed] (c)
            (c) edge (d)
            (d) edge (a)
        ;
        \node at (-1,1/2){$G=$};
        \node at (6.5,1/2){$L_{G,\mathcal{P},\mathcal{Q}} = \left\{\left.\begin{pmatrix}
        a & c & 0 & d \\
        c & a & c & 0 \\
        0 & c & b & d \\
        d & 0 & d & b
        \end{pmatrix} \right\vert a,b,c,d \in \mathbb{R} \right\}$};
        \node at (13,1/2){$\mathcal{M}_{G,\mathcal{P},\mathcal{Q}} = L_{G,\mathcal{P},\mathcal{Q}} \cap \mathcal{S}^n_{++}$};
    \end{tikzpicture}
    \caption{A graph with vertex partition $\mathcal{P} = \{\{1,2\},\{3,4\}\}$, as indicated by the vertex colors,
    and edge partition $\mathcal{Q} = \{\{12,23\},\{34,14\}\}$, as indicated by the edge line styles,
    alongside the corresponding RCON model.
    }\label{fig: RCON example}
\end{figure}
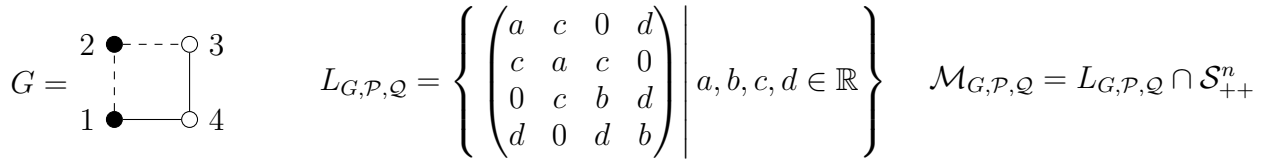

It would be interesting to understand MLTs of non-generic linear concentration models, other than graphical models.
One such family, RCON models, was introduced by H{\o}jsgaard and Lauritzen in~\cite{hojsgaard2008graphical} which we now describe.
Let $G$ be a graph on vertex set $\{1,\dots,n\}$ and edge set $E$.
Let $\mathcal{P}$ and $\mathcal{Q}$ be partitions of $\{1,\dots,n\}$ and $E$.
Define $L_{G,\mathcal{P},\mathcal{Q}} \subseteq \mathcal{S}^n$ to consist of all symmetric $n\times n$ matrices $M$
such that $M_{ij} = 0$ if $ij \notin E$ and $M_{ii} = M_{jj}$ if $i,j$ are in the same part in $\mathcal{P}$
and $M_{ij} = M_{kl}$ for all edge pairs $ij,kl$ in the same part of $\mathcal{Q}$.
Then define
\[
    {\mathcal{M}}_{G,\mathcal{P},\mathcal{Q}} := L_{G,\mathcal{P},\mathcal{Q}} \cap \mathcal{S}^n_{++}.
\]
The \emph{RCON model} corresponding to $(G,\mathcal{P},\mathcal{Q})$ is the set of multivariate Gaussian distributions whose inverse covariance matrix lies in $\mathcal{M}_{G,\mathcal{P},\mathcal{Q}}$.
See figure~\ref{fig: RCON example} for an example.
It is an interesting problem to find methods for computing and bounding the MLTs of RCON models in terms of the combinatorics of the underlying graph and vetex/edge partitions.
Computations of several such examples can be found in~\cite{uhler2012geometry}.

One can also consider Gaussian models defined by things other than linear constraints on the inverse covariance matrix.
A particularly important family of such models are called \emph{Gaussian structural equation models}
which we now describe.
Let $G$ be a directed graph on vertex set $[n]$ and arc set $D$.
One should interpret $G$ as expressing the causal relation among $n$ random variables,
where $i \rightarrow j \in D$ indicates that the random variable corresponding to $i$ has a direct impact on the random variable corresponding to $j$.
The corresponding Gaussian structural equation model $\mathcal{M}_G$ is the model consisting of Gaussian distributions
whose covariance matrices can be expressed as
\[
    ((I-\Lambda)^{-1})^TD (I-\Lambda)^{-1}
\]
where $D$ is a diagonal matrix, $I$ is the identity matrix, $\Lambda_{ij} = 0$ whenever $i\rightarrow j \notin D$.
Thanks to work of Drton, Fox, K\"aufl, and Pouliot~\cite{drton2019maximum}, maximum likelihood thresholds of Gaussian structural equation models are characterized as follows.
The original result was about a slightly more general class of models, but we only discuss this particular case for ease of exposition.

\begin{thm}[{\cite[Theorem 2]{drton2019maximum}}]\label{thm: MLT of directed Gaussian models}
    Let $G$ be a directed graph.
    Then the maximum likelihood threshold of $\mathcal{M}_G$
    is the maximum in-degree of a vertex of $G$.
\end{thm}

\bibliography{graphicalModels}
\bibliographystyle{plain}

\end{document}